\magnification 1200
        \def\R{{\rm I\kern-0.2em R\kern0.2em \kern-0.2em}}
        \def\N{{\rm I\kern-0.2em N\kern0.2em \kern-0.2em}}
        \def\P{{\rm I\kern-0.2em P\kern0.2em \kern-0.2em}}
        \def\B{{\rm I\kern-0.2em B\kern0.2em \kern-0.2em}}
        \def\Z{{\rm I\kern-0.2em Z\kern0.2em \kern-0.2em}}
        \def\C{{\bf \rm C}\kern-.4em {\vrule height1.4ex width.08em depth-.04ex}\;}
        \def\B{{\bf \rm B}\kern-.4em {\vrule height1.4ex width.08em depth-.04ex}\;}
        
        \def\D{{\Delta}}

        \def\z{{\zeta}}

        \def\L{{\cal L}}

        \
        \vskip 20mm
        \centerline {\bf ON HOLOMORPHIC FUNCTIONS WITH CLUSTER SETS OF}

\centerline{\bf FINITE LINEAR MEASURE}
        \vskip 4mm
        \centerline{Josip Globevnik and David Kalaj}
        \vskip 4mm

 \noindent        \bf Abstract\ \ \rm We prove that if  $f$ is a holomorphic function on the open unit disc in $\C$ whose cluster set $C(f)$ has finite linear measure and is such that $\C\setminus C(f)$ has finitely many components, then the derivative $f^\prime$ belongs to the Hardy space  $H^1$. 
        \vskip 4mm
 \noindent        \bf 1. Introduction and the result \rm
\vskip 2mm Let $f$ be a complex valued function defined on the open unit disc $\D\subset\C$. The \it cluster
 set \rm $C(f)$ of $f$ is defined as the set of all $w\in\C$ such that $w=\lim _{n\rightarrow\infty} f(z_n)$ where $z_n\in\D,\ 
\lim _{n\rightarrow\infty} |z_n|=1$. If $f$ is bounded and continuous then $C(f)$ is a compact connected set.

If $E\subset \C$ then  \it the linear measure\rm \ $\Lambda (E)$ of $E$ is the one dimensional Hausdorff measure of $E$, that is, 
$$
\Lambda (E)=\lim _{\varepsilon \rightarrow 0} \inf _{(D_n)} \sum_n \hbox{diam} D_n
$$
 where the infimum is taken over all systems of discs with $\hbox{diam} D_n < \varepsilon$ that cover $E$. If $E$ is
 a rectifiable Jordan arc then $\Lambda (E)$ equals the arclength of $E$.

Let $f$ be a nonconstant holomorphic function on $\Delta$.
\vskip 2mm
\noindent \bf Theorem 1.1\ \rm [D, p.42]\ \it\ The function $f$ extends continuously through $\overline\D$ and $f$ is absolutely continuous on $b\D$ if
 and only if $f^\prime$ belongs to Hardy space $H^1$, that is, if and only if 
$$\sup_{0<r<1}\int_0^{2\pi}|f^\prime (re^{i\theta})|d\theta <\infty .
$$
This happens  if $f$ is of bounded variation on $b\D$. \rm
\vskip 2mm
\noindent A special case when this happens is
\vskip 2mm
\noindent \bf Theorem 1.2\ \rm [P2, p.320]\it \ Let $f$ be a biholomorphic map from $\D$ onto a
 domain $D$. Then $f^\prime\in H^1$ if and only if $\Lambda (bD)<\infty $. In particular, if $D$ is bounded by a rectifiable simple closed curve then $f^\prime\in H^1$.\rm
\vskip 2mm
We want to find more general conditions for a holomorphic function $f$ on $\D$ which imply 
that $f^\prime\in H^1$. Looking at Theorem 1.2  it is a natural question whether for a holomorphic 
function $f$ on $\D$ the assumption $\Lambda (C(f))<\infty$ implies that $f^\prime\in H^1$.  In the case when   $\Lambda (C(f))<\infty$ we have the following result of H.\ Alexander and Ch.\ Pommerenke
\vskip 2mm
\noindent
\bf Theorem 1.3\ \rm [A, P1] \it Let $f$ be a holomorphic function on $\D$ such 
that $\Lambda (C(f))<\infty$. Then $f$ extends continuously through $\overline\D$. \rm
\vskip 2mm
\noindent However, as shown by D.\ Gnuschke-Hauschild [G-H, p.592], the condition that $\Lambda (C(f))<\infty 
$ does not necessarily imply that $f^\prime \in H^1$. 

If $f$ is a nonconstant holomorphic function on $\D$ then $f(\D ) $ is a nonempty open set; if $\Lambda (C(f))<\infty $ then $C(f)$ is a 
compact connected set with empty interior so $\C\setminus C(f)$ is an open set having at least one bounded component. 
 In the example of Gnuschke-Hauschild $\C\setminus C(f)$ has infinitely many components. The pricipal result of the present paper
 is that the situation is different when  $\C\setminus C(f)$ has only finitely many components:
\vskip 2mm
\noindent \bf Theorem 1.4\ \it Let $f$ be a holomorphic function on $\D $ such that  $\Lambda (C(f))<\infty 
$ and such that $\C\setminus C(f)$ has finitely many components. Then $f^\prime\in H^1$.
\vskip 4mm
\noindent \bf 2. Proof of Theorem 1.4, Part 1 \rm
\vskip 2mm Let $f$ be as in Theorem 1.4. With no loss of generality assume that $f$ is not a constant. 
By Theorem 1.3, f extends continuously through $\overline\D$. 
We denote the extension with the same letter $f$. By Theorem 1.1 we
 have to show that there is a constant $M<\infty$ such that
$$
\sum_{j=1}^\ell |f(e^{i\theta _j})- f(e^{i\theta_{j-1}})|< M
$$
whenever $0\leq \theta_0<\theta_1<\cdots <\theta_\ell\leq 2\pi$.

As Gnuschke-Hauschild did in  [G-H] we decompose $\D$ into subsets on which $f$ has simpler behavior.
 If $V$ is a component ot the open set $f(\D )\setminus C(f)$ then V is an open connected set whose boundary 
is contained in $C(f)=f(b\D )$ and hence it is a component of $\C\setminus C(f)$.  By our assumption  $\C\setminus C(f)$ has 
finitely many components so it follows that $f(\D)\setminus C(f)$ has finitely many components. Denote these components by
 $U_1, U_2,\cdots U_m$. Since $C(f)$ is connected the domains $U_j$ are simply connected.

The set $E= (f|\D )^{-1}(C(f)) = \{ z\in\D\colon f(z)\in C(f)\} = \{ z\in\D\colon f(z)\in f(b\D )\}$ is a closed subset of $\D$. Let $G_k, \ k=1,2,\cdots$,\ 
be the components of the open set $\D\setminus E$. Then, as shown in [G-H], each domain $G_k$ is simply connected and $f(G_k)$ equals one of the 
domains  $U_1,\cdots , U_m$, moreover $f|G_k\colon G_k\rightarrow f(G_k)$ is a  proper map, that is, if $\varphi_k\colon\D\rightarrow G_k$ and  $\psi_k
\colon \D\rightarrow f(G_k)$ are biholomorphic maps then $g_k=\psi_k^{-1}\circ f\circ \varphi_k$ is a proper holomorphic map from $\D$ to 
$\D$, that is, a finite Blaschke product whose multiplicity we denote by $\nu_k$.

We show that the number of components $G_k$ is finite. To see this, assume the opposite, that there are one of the components $U$ among $U_1,\cdots , U_m$ and a sequence of components $G_k, \ k\in\N$, such that for each $k\in\N,\ f(G_k)=U$. Let $w\in U$. Clearly $w\not\in C(f)$. For each $k$ there is a $z_k\in G_k$ such that $f(z_k) = w$. Since $G_k$ are pairwise disjoint, $z_k$ is an injective sequence. We  claim that $|z_k|\rightarrow 1$ as $k\rightarrow\infty$. If not, passing to a subsequence if necessary, we may assume that $z_k\rightarrow z\in\D$ and thus $f(z_k)=w$ for all k which is not possible since $f$ is nonconstant. Thus, $|z_k|\rightarrow 1$. But since $f(z_k)=w$ for all $k$ it follows that  $w\in C(f)$, a contradiction. Thus, the number of components $G_k$ 
is finite, denote them by $G_1, G_2,\cdots G_n$. Obviously $n\geq m$.

The cluster set $C(f)$ is a plane continuum. Since $\Lambda (C(f))<\infty$ it follows that $C(f)$ is locally connected [CC, Lemma 2, p.49]. Since  $\C\setminus 
 C(f)$ has finitely many components, Thorhorst's theorem [CC, p.44, \ W, p.113]  implies that each component of $\C\setminus C(f)$ has locally connected boundary. Thus, $bU_j$ is locally connected for each $j,\ 1\leq j\leq m$.

For each $G_j$ there is some $U_k$ such that $f(G_j)=U_k$ and such that $f|G_j\colon G_j\rightarrow U_k$ is a proper map. Since  $bU_k$ is locally connected, [CC, Lemma 1, p.46] implies that $bG_j$ is locally connected. Thus, for each $j,\ 1\leq j\leq n$, the domain $G_j$ is simply connected and has locally connected boundary. 

For each $j,\ 1\leq j\leq n$, let $\psi_j\colon\D\rightarrow G_j$ be a biholomorphic map. Since $bG_j$ is locally connected, the map $\psi_j$ extends continuously through $\overline\D$ and $\psi_j(b\D)=bG_j$ \ [P2, p.\ 279]. In particular, for each $a,b\in bG_j,\ a\not=b$, there is  an arc $L$, contained in $G_j$ except for its endpoints $a$ and $b$. 

Notice that the set $E=(f|\D )^{-1}(C(f))$ has no interior: If $E$ contains a nonempty open set then, since $f$ is open, $f(E)$ contains a nonempty open set. Since $f(E)\subset C(f)$ this is not possible since $C(f)$, being of finite linear measure, has no interior. Since $E$ has no interior and since $E\cup [\cup_{j=1}^n G_j]=\D$ it follows that $\cup_{j=1}^n G_j$ is dense in $\D$. Thus, $\cup_{j=1}^n\overline{G_j}=\overline\D$ so $b\D\subset \cup _{j=1}^n bG_j.$
\vskip 4mm
\noindent \bf 3. Some estimates
\vskip 2mm \rm
Assume that $D$ is a domain such that  $\Lambda (bD)<\infty$ and that $\Phi\colon\D\rightarrow D$ is a biholomorphic map. We know that $\Phi$ extends continuously through $\overline\D$. The proof of [Pom., Th.10.1, p.321] shows that 
$$ 
\int_0^{2\pi}|\Phi^\prime (re^{i\theta})|d\theta \leq \pi\Lambda (bD)\ \ (0<r<1).
\eqno (3.1)
$$
Since\ \   $r\int_0^{2\pi}|\Phi^\prime (re^{i\theta})|d\theta $ \ equals the length of the curve 
$\{ \theta\mapsto\Phi (re^{i\theta})\colon 0\leq \theta\leq 2\pi\} $ it follows that given 
$\theta_j,$
$$
\theta_0<\theta_1<\cdots <\theta _\ell \leq \theta_0+2\pi
\eqno (3.2)
$$
we have $\sum_{j=1}^\ell |\Phi(re^{i\theta _j}) - \Phi(re^{i\theta _{j-1}})|\leq 
r\int_0^{2\pi}|\Phi^\prime (re^{i\theta})| d\theta \leq r\pi \Lambda (bD)$,
which, as $r\rightarrow 1$, gives
$$
\sum_{j=1}^\ell |\Phi(e^{i\theta _j})-\Phi (e^{i\theta_{j-1}})| \leq \pi\Lambda (bD)
\eqno (3.3)
$$
whenever $\theta_j$ satisfy (3.2).

Assume now that $D$ is a simply connected domain with $\Lambda (bD)<\infty$ and that 
$g\colon \D\rightarrow D$ is a proper holomorphic map. If $\Phi\colon \D\rightarrow D$ is a biholomorphic map then $\Phi^{-1}\circ g$ is a proper holomorphic map from $\D$ to $\D$, a finite Blaschke product $B$ , so $g=\Phi\circ B$ where the multiplicity of the map $g$ equals the multiplicity of $B$.

If $\nu$ is the multiplicity of $B$ then, as $\z $ runs around $b\D$ once, $B(\z)$ runs around
$b\D$ $\nu$  times so if $\theta_j$ satisfy (3.2) then there are $\tau_j,\ \tau_0<\tau_1<\cdots< \tau_\ell \leq \tau_0+\nu 2\pi$ such that $B(e^{i\theta_j}) = e^{i\tau_j}
\ (1\leq j\leq \ell)$. The preceding discussion now implies that
$\sum_{j=1}^{\ell} |g(e^{i\theta _j})-g(^{i\theta_{j-1}})| = \sum_{j=1}^\ell|\Phi(e^{i\tau_j})-\Phi(e^{i\tau_{j-1}})| \leq \nu\pi\Lambda(bD)$. This gives
\vskip 2mm
\noindent\bf Proposition 3.1 \it Let $D$ be a simply connected domain such that $\Lambda (bD)<\infty$ and let  $g\colon \D\rightarrow D$ be a proper holomorphic map of multiplicity $\nu$. Then $g$ extends continuously through $\overline\D$. If $\theta_0<\theta_1<\cdots <\theta_\ell \leq \theta_0+2\pi$ then
$$
\sum_{j=1}^\ell |g(e^{i\theta _j})-g(e^{i\theta_{j-1}})|\leq \nu\pi\Lambda (bD)
$$
\vskip 4mm
\noindent\bf 4. Proof of Theorem 1.4, Part 2 \rm
\vskip 2mm
We shall show that
$$\left. \eqalign{
&b\D \hbox{\ can be written as a finite union of pairwise disjoint semiopen arcs\ }\cr
& L_1, L_2,\cdots ,L_\mu, \hbox{\  each being of the form\ }
 \{ e^{i\theta}\colon \alpha\leq \theta<\beta\}\hbox{\  such that each\ }\cr
 &L_k \ \hbox{ \   is contained in\ } bG_\sigma \hbox{\ for some \ }\sigma,\ 1\leq\sigma\leq n. 
\cr}
\right\}
\eqno (4.1)
$$
Asssume  for a moment that we have done this. With no loss of generality assume that the initial point $e^{i\alpha}$ of $L_1$ is $1$. 

Let $A_j,\ 1\leq j\leq \nu$, be points on $b\D$. We say that the $\nu-$tuple $(A_1,A_2,\cdots, A_\nu)$ is \it ordered positively \rm if $A_j=e^{i\theta_j}$ where $ \theta_1<\theta_2<\cdots <\theta_\nu <\theta_1+2\pi$.

To
proceed we need the following
\vskip 2mm\noindent\bf Proposition 4.1\ \it Let $G\subset\D$ be a simply connected domain with locally connected boundary and let  $\Phi\colon\D \rightarrow G$ be 
a biholomorphic map (that extends continuously through $\overline\D$ as $bG$ is locally connected). Suppose that $L$ is a closed arc in $b\D$ such that $L\subset bG$ and assume that $A_j\in L,\ 1\leq j \leq\nu$, are such that the $\nu-$tuple $(A_1,A_2,\cdots, A_\nu)$ is \it ordered positively. 
For each $j,\ 1\leq j\leq\nu,$ let $a_j\in b\D $ be such that $\Phi (a_j)=A_j\ (1\leq j\leq \nu)$. Then the $\nu-$tuple $(a_1,\ a_2,\cdots , a_\nu)$  is ordered positively.\rm
\vskip 2mm 
\noindent\bf Remark.\ \rm Note that we do not assume that $\hbox{Int}L$ is an open subset of $bG$. Note also that since $\Phi (b\D)=bG$ it follows that for each $j,\ 1\leq j\leq \nu,$ \ there is an $a_j\in b\D$ which satisfies $\Phi(a_j)=A_j$. Note that this $a_j$ is not necessarily unique. 
\vskip 2mm
\noindent \bf Proof. \rm It is enough to prove the special case of Proposition 4.1  when $\nu=3$ as the general case will then follow by using this special case inductively. So let $\nu=3$ and let $A_j,\ a_j,\ 1\leq j\leq 3$, be as in Proposition 4.1. Let $\lambda\subset b\D$  be the arc obtained by sliding $a_1$ to $a_3$ along $b\D$ in positive direction, that is, if $a_1=e^{i\omega_1}$ and $a_3=e^{i\omega_3}$ where $\omega_1<\omega_3<\omega_1+2\pi$, then $
\lambda=\{e^{i\omega}\colon\ \omega_1<\omega<\omega_3\}$. To see  that $(a_1,a_2,a_3)$ is ordered positively we must show that $a_2\in\lambda$.

Let $\ell$ be the arc in $\D$ consisting of the segment joining $a_3$ with $0$ and the segment joining $0$ with $a_1$ and let $\Omega $ be the domain bounded by $\lambda\cup\ell$. Orient $b\Omega =\lambda\cup \ell$ in the positive direction. In particular, $\ell $ has $a_3$ as the initial point and $a_1$ as the final point. The arc $\L=\Phi (\ell )$ is contained in $\D$ except for its endpoints $A_3=\Phi(a_3)$ and $A_1=\Phi(a_1)$; we keep the orientation from $\ell $ so $A_3$ is the initial point and $A_1$ is the final point of $\L$.

Since the endpoints of $\L$ belong to $b\D$,\ $\D\setminus \L$ has two components $D_1$ and $D_2$ where $D_1$ is bounded by $\L$ and by the arc $\Sigma\subset b\D$ obtained  by sliding $A_1$ to $A_3$ along $b\D$ in positive direction. The arc $\Sigma$ oriented in this direction together with the arc $\L$ oriented as above form the positively oriented boundary of $D_1$. Since $\Phi\colon\D\rightarrow \Phi(\D) =G$ is a biholomorphic  map, $\Phi(\Omega)$ must be contained  either in $D_1$ or in $D_2$. Along $\ell,\ \Omega$ lies to  the left of $\ell$ with respect to the above orientation of $\ell$ and since $\Phi$ is conformal it follows that  along $\L$, $\Phi(\Omega)$ lies to the left of $L$. This implies that $\Phi (\Omega)\subset D_1$.  If $a_2 \in b\D\setminus\lambda$ then $a_2$ is in the closure of $\D\setminus\Omega$, which, by the continuity of $\Phi$ on $\overline\D$, implies that $\Phi (a_2) = A_2$ is in the closure of $D_2$ which is impossible since 
$A_2$ is an interior point of $\Sigma$. This completes the proof of Proposition 4.1.
\vskip 2mm
Recall that for each $\sigma,\ 1\leq\sigma\leq n$, \ $\nu_\sigma$ is the multiplicity of the proper map $f|G_\sigma\colon G_\sigma \rightarrow f(G_\sigma)$.  Let $V=\max_{1\leq \sigma\leq n} \nu_\sigma$. 

Now, let
$$
0\leq \theta_0<\theta_1<\cdots<\theta_\ell\leq 2\pi.
\eqno (4.2))
$$
Fix $k,\ 1\leq k\leq \mu$, and consider those points $e^{i\theta_j}$ that belong to $L_k$, suppose that these points are $e^{i\theta_s}, e^{i\theta_{s+1}},\cdots, e^{i\theta_t}$. Let
$\sigma,\ 1\leq\sigma\leq n$, be such that $L_k\subset bG_\sigma$. By Proposition 4.1 there are $\omega_j,\ \omega_s<\omega_{s+1}<\cdots <\omega_t <\omega _s+2\pi$, such that 
$$
\psi_\sigma (e^{i\omega_p})=e^{i\theta_p}\ \ (s\leq p\leq t)
$$
and hence by Proposition 3.1
$$
\eqalign
{
\sum_{j=s+1}^t |f(e^{i\theta_j})-f(e^{i\theta_{j-1}})| =
\sum_{j=s+1}^t |f(\psi_\sigma(e^{i\omega_j}))-f(\psi_\sigma(e^{i\omega_{j-1}}))| \leq  \cr
\leq\nu_\sigma\pi\Lambda (b(f(G_\sigma))) 
\leq V \pi\Lambda (C(f))  .\cr}
$$
We repeat this for every $k, \ 1\leq k\leq \mu$.  In the sum
$$
\sum _{j=1}^\ell |f(e^{i\theta_j})-f(e^{i\theta_{j-1}})|
\eqno (4.3)
$$
there are terms as above whose total sum does not exceed $\mu V\pi\Lambda (C(f))$. However, in the sum (4.3) there are also terms  $|f(e^{i\theta_j})-f(e^{i\theta_{j-1}})|$  where $e^{i\theta_{j-1}}$ is the last point in $L_k$ and $e^{i\theta_j}$ is the first point in 
$L_{k+1}$. Each such term does not exceed $\hbox{diam}(C(f))$ and there are $\mu -1$ such terms. Thus, the sum (4.3) is bounded above by $\mu V\pi(\Lambda (C(f)) + (\mu-1)\hbox{diam}C(f)$ whenever $\theta_j$ satisfy (4.2). Consequently, $f|(b\D)$ has
bounded variation which was to be proved. 

It remains to prove (4.1).
\vskip 4mm
\noindent\bf 5. Proof of Theorem 1.4, Part 3 \rm
\vskip 2mm
We now want to study how domains $G_j$ "touch" $b\D$, that is we want to look at $(b\D)\cap \overline{G_j} = (b\D)\cap bG_j$. Fix $j, \ 1\leq j\leq n$.  Consider the set $\Sigma = \overline{G_j}\cap b\D = (bG_j)\cap b\D$. This is a closed subset of $b\D$. Its complement 
$(b\D)\setminus \Sigma$ is an open subset of $b\D$. We show that it has finitely many components. To see this, let $\lambda\subset b\D$ be a component of $b\D\setminus\Sigma$ Thus, $\lambda $ is an open arc on $b\D$ disjoint from $\overline{G_j}$ whose endpoints belong to $\overline{G_j}$. As observed in Section 2 there is an arc L contained in $G_j$ except for its endpoints which coincide with endpoints of $\lambda$. Since $b\D=(\cup_{i=1}^nbG_i)\cap b\D$ it follows that there are a $k,\ 1\leq k\leq n,\ k\not=j$, and a point $a\in\lambda\cap (bG_k)$. Notice that $bG_k$ meets no other component of $(b\D)\setminus\Sigma$. Indeed, if it does, there is an arc $L_1$ contained in $G_k$ except its 
endpoints $a$ and $b\in\lambda_1$ where $\lambda_1$ is a component of $b\D\setminus\Sigma$ different from $\lambda$. However, such an $L_1$ would have to intersect $L$ which is impossible since $L\subset G_j$ and $L_1\subset G_k$ and $G_j\cap G_k = \emptyset$. Thus, for each component $\lambda$ of $b\D\setminus\Sigma$ there is a $G_k$ such that $\lambda\cap bG_k\not=\emptyset $ and such that $bG_k$ meets no other 
component of $b\D\setminus\Sigma$. Obviously all domains $G_k$ so obtained are pairwise disjoint. Since there are only $n$ of them it follows that $b\D\setminus\Sigma$ has finitely many components. 

Thus, for each $j,\ 1\leq j\leq n$, the set $bG_j\cap b\D$, if nonempty is a union of finitely many pairwise disjoint closed sets each of which is either a point or a closed arc. (The only trivial noninteresting exception is when there is just one such $G_j$, \ $G_1=\D$ when $bG_1=b\D$.) Since $\cup_{j=1}^n(bG_j\cap b\D)=b\D$ it follows that the family of all closed arcs so obtained for all $j,\ 1\leq j\leq n$, covers $b\D$. It is easy to see that this implies (4.1). The proof of Theorem 1.4 is complete.
\vskip 6mm
\centerline{\bf References}
\vskip 3mm
\noindent [A] H.\ Alexander: Polynomial hulls and linear measure.  Complex analysis, II (College Park, Md., 1985–86),  1-11, Lecture Notes in Math., 1276, Springer, Berlin, 1987. 
\vskip 2mm
\noindent [CC] J.\ J.\ Carmona and J.\ Cufi:  On analytic functions with locally connected cluster sets. Complex Variables Theory Appl.  10  (1988),  no. 1, 43-52.
\vskip 2mm
\noindent [G-H] D.\  Gnuschke-Hauschild: On the angular derivative of analytic functions. Math. Z.  196  (1987),  no. 4, 591-601.
\vskip 2mm
\noindent [D] P.\ L.\ Duren: \it Theory of $H^p$   spaces.\rm\  Pure and Applied Mathematics, Vol. 38 Academic Press, New York-London 1970
\vskip 2mm
\noindent 
[P1] Ch.\ Pommerenke:  On analytic functions with cluster sets of finite linear measure. Michigan Math. J.  34  (1987),  no. 1, 93-97
\vskip 2mm
\noindent 
[P2] Ch.\ Pommerenke: \it Univalent functions. With a chapter on quadratic differentials by Gerd Jensen.\rm \ Studia Mathematica/Mathematische Lehrbücher, Band XXV. Vandenhoeck \& Ruprecht, Goettingen, 1975
\vskip 2mm
\noindent
[W] G.\ T.\ Whyburn: \it Analytic topology.\rm \  American Mathematical Society Colloquium Publications, Vol. XXVIII American Mathematical Society, Providence, R.I. 1963 
\vskip 8mm
\noindent Institute of Mathematics, Physics and Mechanics and

\noindent  Department of Mathematics, University of Ljubljana

\noindent Jadranska 19, 1000 Ljubljana, Slovenia  (Globevnik)
\vskip 3mm
\noindent Faculty of Natural Sciences and Mathematics, University of Montenegro

\noindent Cetinjski put b.b., 81000 Podgorica, Montenegro (Kalaj)

\

\bye